\documentclass[11pt]{article}
\usepackage{indentfirst}
\usepackage{amsfonts}
\usepackage{url}
\usepackage{algorithm}
\usepackage{algpseudocode}
\usepackage{listings}
\newtheorem{teor}{Theorem}[section]
\newtheorem{exemplu}{Example}[section]
\author{E. Scheiber\thanks{retired from \textit{Transilvania} University of Bra\c{s}ov, e-mail: scheiber@unitbv.ro}}
\title{On the Convergence of the Variational Iteration Method}
\date{}
\begin{document}
\maketitle

\begin{abstract}
Convergence results are stated for the variational iteration method applied to
solve an initial value problem for a system of ordinary differential equations.


\end{abstract}

\section{Introduction}
The Ji-Huan He's Variational Iteration Method (VIM) was applied to a large range of problems for ordinary and partial 
differential equations. The purpose of this paper is to prove convergence theorem 
when the method is applied to solve an initial value problem for a system of ordinary differential
equations. In advance there is presented a convergence result in the case of an initial value
problem for an ordinary differential equation. The result as well the proof may be found in 
Salkuyeh D. K., Tavakoli A., \cite{5}, too. In the case of a system of linear differential equations
a convergence result is given by Salkuyeh D. K., \cite{3}.
Variants of the convergence of the VIM are studied in \cite{4}, \cite{6}, \cite{7}.

In the last section there are presented some results of our computational experiences.

\section{The case of an ordinary differential equation}

Let $f:[t_0,t_f]\times\mathbb{R}\rightarrow \mathbb{R}\ (t_0<t_f<\infty),$ be a continuous  function such that it has
the first and second order partial derivatives in the second variable. Denoting 
$\frac{\partial f(t,x)}{\partial x}=f_x(t,x),$ 
we suppose that there exists $L>$ such that
$$
|f_x(t,x)|\le L \qquad \forall (t,x)\in [t_0,t_f]\times \mathbb{R}.
$$
As a consequence the function $f$ satisfies the Lipschitz property 
$$
|f(t,x)-f(t,y)|\le L|x-y|,\qquad \forall x,y\in\mathbb{R}.
$$

To solve the initial value problem
\begin{eqnarray}
x'(t)&=&f(t,x(t)) \quad  t\in[t_0,t_f] \label{vim1}\\
x(t_0)&=&x^0 \label{vim2}
\end{eqnarray}
the VIM is applied, \cite{2}.
There is considered the sequence 
\begin{equation}\label{vim3}
u_{n+1}(t)=u_n(t)+\int_{t_0}^t\lambda(s)(u'_n(s)-f(s,u_n(s)))\mathrm{d}s, \quad n\in\mathbb{N}.
\end{equation}
where $\lambda$ is the so called the Lagrange multiplier \cite{1}.
There is supposed that $u_n(t_0)=x^0$ and that $u_n$ is a continuous derivable function.

Denoting  $x(t)$ the solution of the initial  value problem (\ref{vim1})-(\ref{vim2}), if
$u_n(t)=x(t)+\delta u_n(t)$ and $u_{n+1}(t)=x(t)+\delta u_{n+1}(t)$ then (\ref{vim3}) implies
$$
\delta u_{n+1}(t)=\delta u_n(t)+\int_{t_0}^t\lambda(s)\left(x'(s)+\delta u_n'(s)-f(s,x(s)+\delta u_n(s)\right)\mathrm{d}s=
$$
$$
=\delta u_n(t)+\int_{t_0}^t\lambda(s)\left(x'(s)+\delta u_n'(s)-f(s,x(s))-f_x(s,x(s))\delta u_n(s)\right)\mathrm{d}s+O((\delta u_n)^2)=
$$
$$
=\delta u_n(t)+\int_{t_0}^t\lambda(s)\left(\delta u_n'(s)-f_x(s,x(s))\delta u_n(s)\right)\mathrm{d}s+O((\delta u_n)^2).
$$
After an integration by parts the above equality becomes
$$
\delta u_{n+1}(t)=(1+\lambda(t))\delta u_n(t)-\int_{t_0}^t\left(\lambda'(s)+f_x(s,x(s))\lambda(s)\right)\delta u_n(s)\mathrm{d}s+O((\delta u_n)^2).
$$
In order that $u_{n+1}$ be a better approximation than $u_n,$ there is required that $\lambda$ to be the solution
of the following initial value problem
\begin{eqnarray}
\lambda'(s)&=&-f_x(s,x(s))\lambda(s)\quad s\in[t_0,t] \label{vim4}\\
\lambda(t)&=&-1\label{vim5}
\end{eqnarray}
Because $x(s)$ is an unknown function, instead of (\ref{vim4})-(\ref{vim5}) there is 
considered the problem
\begin{eqnarray}
\lambda'(s)&=&-f_x(s,u_n(s))\lambda(s)\quad s\in[t_0,t] \label{vim6}\\
\lambda(t)&=&-1\label{vim7}
\end{eqnarray}
with the solution denoted  $\lambda_n(s,t).$

The recurrence formula (\ref{vim3}) becomes
\begin{equation}\label{vim8}
u_{n+1}(t)=u_n(t)+\int_{t_0}^t\lambda_n(s,t)(u'_n(s)-f(s,u_n(s)))\mathrm{d}s, \quad n\in\mathbb{N}.
\end{equation}

The solution of the initial value problem  (\ref{vim6})-(\ref{vim7}) is 
\begin{equation}\label{vim9}
\lambda_n(s,t)=-e^{\int_s^tf_x(\tau,u_n(\tau))\mathrm{d}\tau}.
\end{equation}
It results that $|\lambda_n(s,t)|\le e^{L(t-s)}\le e^{LT},\ \forall\ t_0\le s\le t\le t_f$ and $T=t_f-t_0.$

The following convergence result occurs:
\begin{teor}\label{vim23}\cite{5}(Th. 2)
If $f:[t_0,t_f]\times\mathbb{R}\rightarrow \mathbb{R}$ is a continuous  function such that it has
the first order partial derivatives in $x$ bounded  $|f_x(t,x)|\le L, \forall (t,x)\in[t_0,t_f]\times\mathbb{R}$
then the sequence $(u_n)_{n\in\mathbb{N}}$ defined by (\ref{vim8}) converge uniformly to $x(t),$ the
solution of the initial value problem  (\ref{vim1})-(\ref{vim2}).
\end{teor}

\vspace{0.3cm}\noindent
\textbf{Proof.}
Subtracting the equality
$$
x(t)=x(t)+\int_0^t\lambda_n(s,t)\left(x'(s)-f(s,x(s))\right)\mathrm{d}s
$$
from (\ref{vim8}) it results
$$
e_{n+1}(t)=e_n(t)+\int_{t_0}^t\lambda_n(s,t)\left(e'_n(s)-(f(s,u_n(s))-f(s,x(s)))\right)\mathrm{d}s,
$$
where $e_n(t)=u_n(t)-x(t),\ n\in\mathbb{N}.$

Again, after an integration by parts, it results
$$
e_{n+1}(t)=-\int_{t_0}^t\left(\lambda'_n(s,t)e_n(s)+(f(s,u_n(s))-f(s,x(s)))\lambda_n(s,t)\right)\mathrm{d}s.
$$
Taking into account (\ref{vim6}), the above equality gives
$$  
e_{n+1}(t)=\int_{t_0}^t\lambda_n(s,t)\left(f_x(s,u_n(s))e_n(s)-(f(s,u_n(s))-f(s,x(s)))\right)\mathrm{d}s.
$$  
The hypothesis on $f$ implies the inequality
$$
|f_x(s,u_n(s))e_n(s)-(f(s,u_n(s))-f(s,x(s))|\le 2L|e_n(s)|
$$
end consequently
\begin{equation}\label{vim10}
|e_{n+1}(t)|\le 2Le^{LT}\int_{t_0}^t|e_n(s)|\mathrm{d}s.
\end{equation}

Let be $M=2Le^{LT}.$ For any continuous function $\varphi$ on $[t_0,t_f]$ we use the notation
$\|\varphi\|_{\infty}=\max_{t_0\le t\le t_f}|\varphi(t)|.$

From (\ref{vim10}) we obtain successively:

\noindent
For $n=0$ 
$$
|e_1(t)|\le M \int_{t_0}^t|e_0(s)|\mathrm{d}s\le M(t-t_0)\|e_0\|_{\infty} \ \Rightarrow\ \|e_1\|_{\infty}\le MT\|e_0\|_{\infty}.
$$

\noindent
For $n=1$ 
$$
|e_2(t)|\le M \int_{t_0}^t|e_1(s)|\mathrm{d}s\le \frac{M^2 (t-t_0)^2}{2}\|e_0\|_{\infty}\ \Rightarrow\ \|e_2\|_{\infty}\le \frac{M^2T^2}{2}
\|e_0\|_{\infty}.
$$

\noindent
Inductively, it results that
$$
|e_n(t)|\le M\int_{t_0}^t|e_{n-1}(s)|\mathrm{d}s\le \frac{M^n (t-t_0)^n}{n!}\|e_0\|_{\infty}\ \Rightarrow\ \|e_n\|_{\infty}\le \frac{M^nT^n}{n!}
\|e_0\|_{\infty}.
$$
and hence $\lim_{n\rightarrow\infty}\|e_n\|_{\infty}=0.\quad\rule{5pt}{5pt}$

\section{The case of a system of ordinary differential equations}

Let be the system of ordinary differential equation 
\begin{equation}\label{vim20}
\left\{\begin{array}{lcl}
x_1'(t)=f_1(t,x_1(t),\ldots,x_m(t)) & \quad & x_1(t_0)=x_1^0\\
\vdots \\
x_m'(t)=f_m(t,x_1(t),\ldots,x_m(t)) & \quad & x_m(t_0)=x_m^0
\end{array}\right.
\end{equation}
where $t\in [t_0,t_f]$ with $t_0<t_f<\infty.$

We shall use the notations
\begin{eqnarray*}
\mathbf{x}=(x_1,\ldots,x_m)\\
\|\mathbf{x}\|_1=\sum_{j=1}^m|x_j| \\
\|\mathbf{x}\|_{\infty}=\max_t\|\mathbf{x}(t)\|_1
\end{eqnarray*}

Thus, any equation of (\ref{vim20}) may be rewritten as
$$
x_i'(t)=f_i(t,\mathbf{x}(t)),\quad i\in\{1,\ldots,m\}.
$$

The following hypothesis are introduced: 
\begin{itemize}
\item 
The functions $f_1,\ldots,f_m$ are continuous and have first and second order partial derivatives
in $x_1,\ldots,x_m.$
\item
There exists $L>0$ such that for any $i\in\{1,\ldots,m\}$
$$
|f_i(t,\mathbf{x})-f_i(t,\mathbf{y})|\le L\sum_{j=1}^m|x_i-y_i|=L\|\mathbf{x}-\mathbf{y}\|_1, 
\quad \forall\ \mathbf{x},\mathbf{y}\in\mathbb{R}^m.
$$
As a consequence
$$
|\frac{\partial f_i(t,\mathbf{x})}{\partial x_j}|=|f_{i_{x_j}}(t,\mathbf{x})|\le L,
\quad \forall (t,\mathbf{x})\in [t_0,t_f]\times\mathbb{R}^m,\ \forall i,j\in\{1,\ldots,m\}.
$$
\end{itemize}

According to the VIM there are considered the sequences 
\begin{equation}\label{vim21}
u_{n+1,i}(t)=u_{n,i}(t)+\int_{t_0}^t\lambda_i(s)(u'_{n,i}(s)-f_i(s,\mathbf{u}_n(s)))\mathrm{d}s, \quad n\in\mathbb{N},
\end{equation}
$i\in\{1,\ldots,m\}$ and where $\mathbf{u}_n=(u_{n,1},\ldots,u_{n,m})^T.$ 

There are supposed that $u_{n,i}(t_0)=x^0_i$ and that $u_{n,i}$ is a continuous derivable function for any $i\in\{1,\ldots,m\}.$

The VIM in this case is  little tricky, \cite{3}. The Lagrange multiplier attached to the $i$-th
equation will act only on $x_i.$ 

Denoting  $\mathbf{x}(t)$ the solution of the initial  value problem (\ref{vim20}), if
$u_{n,i}(t)=x_i(t)+\delta u_{n,i}(t)$ and $u_{n+1,i}(t)=x_i(t)+\delta u_{n+1,i}(t)$ but for $j\not=i,\ u_{n,j}(t)=x_j(t)$
 then (\ref{vim21}) implies
$$
\delta u_{n+1,i}(t)=\delta u_{n,i}(t)+\int_{t_0}^t\lambda_i(s)\left(x_i'(s)+\delta u_{n,i}'(s)-\right.
$$
$$
\left.-f_i(s,x_1(s),\ldots,x_{i-1}(s),x_i(s)+\delta u_{n,i}(s),x_{i+1}(s),\ldots,x_m(s))\right)\mathrm{d}s=
$$
$$
=\delta u_{n,i}(t)+\int_{t_0}^t\lambda_i(s)\left(x_i'(s)+\delta u_{n,i}'(s)-f_i(s,\mathbf{x}(s))-f_{i_{x_i}}(s,\mathbf{x}(s))\delta u_{n,i}(s)\right)\mathrm{d}s+
$$
$$
+O((\delta u_{n,i})^2)=
$$
$$
=\delta u_{n,i}(t)+\int_{t_0}^t\lambda_i(s)\left(\delta u_{n,i}'(s)-f_{i_{x_i}}(s,\mathbf{x}(s))\delta u_{n,i}(s)\right)\mathrm{d}s+O((\delta u_{n,i})^2).
$$
Proceeding as in the previous section we find  $\lambda_i(s):=\lambda_{n,i}(s,t)$ as the solution of the initial value problem
\begin{eqnarray*}
\lambda'(s)&=&-f_{i_{x_i}}(s,\mathbf{u}_n(s))\lambda(s)\quad s\in[t_0,t] \\
\lambda(t)&=&-1
\end{eqnarray*}
Then
$$
\lambda_{n,i}(s,t)=-e^{\int_s^tf_{i_{x_i}}(\tau,\mathbf{u}_n(\tau))\mathrm{d}\tau}.
$$
and
$$
|\lambda_{n,i}(s,t)|\le e^{L(t-s)}\le e^{LT},\ \forall\ t_0\le s\le t\le t_f\ \mbox{and}\ T=t_f-t_0.
$$

The recurrence formula (\ref{vim21}) becomes
\begin{equation}\label{vim22}
u_{n+1,i}(t)=u_{n,i}(t)+\int_{t_0}^t\lambda_{n,i}(s,t)(u'_{n,i}(s)-f_i(s,\mathbf{u}_n(s)))\mathrm{d}s, \quad n\in\mathbb{N},
\end{equation}
for any $i\in\{1,\ldots,m\}.$

The convergence result is:
\begin{teor}
If the hypothesis stated above are valid then
the sequence $(\mathbf{u}_n)_{n\in\mathbb{N}}$ defined by (\ref{vim22}) converge uniformly to $\mathbf{x}(t),$ the
solution of the initial value problem  (\ref{vim20}).
\end{teor}

\vspace{0.3cm}\noindent
\textbf{Proof.}
The proof is similar to the proof of Theorem \ref{vim23}.
Subtracting the equality
$$
x_i(t)=x_i(t)+\int_0^t\lambda_{n,i}(s,t)\left(x_i'(s)-f_i(s,\mathbf{x}(s))\right)\mathrm{d}s
$$
from (\ref{vim22}) it results
$$
e_{n+1,i}(t)=e_{n,i}(t)+\int_{t_0}^t\lambda_{n,i}(s,t)\left(e'_{n,i}(s)-(f_i(s,\mathbf{u}_n(s))-f_i(s,\mathbf{x}(s)))\right)\mathrm{d}s,
$$
where $\mathbf{e}_n(t)=\mathbf{u}_n(t)-\mathbf{x}(t),\ n\in\mathbb{N}.$

Again, after an integration by parts, it results
$$  
e_{n+1,i}(t)=\int_{t_0}^t\lambda_{n,i}(s,t)\left(f_{i_{x_i}}(s,\mathbf{u}_n(s))e_{n,i}(s)-(f_i(s,\mathbf{u}_n(s))-f_i(s,\mathbf{x}(s)))\right)\mathrm{d}s.
$$  
The hypothesis on $f_i$ implies the inequality
$$
|f_{i_{x_i}}(s,\mathbf{u}_n(s))e_{n,i}(s)-(f_i(s,\mathbf{u}_n(s))-f_i(s,\mathbf{x}(s))|\le 
 L|e_{n,i}(s)|+L\|\mathbf{e}_n(s)\|_1
 $$
end consequently
$$
|e_{n+1,i}(t)|\le Le^{LT}\int_{t_0}^t(|e_{n,i}(s)|+\|\mathbf{e}_n(s)\|_1)\mathrm{d}s.
$$
Summing these inequalities, for $i=1:m,$ it results
 \begin{equation}\label{vim24}
\|\mathbf{e}_{n+1}(t)\|_1\le (m+1)Le^{LT}\int_{t_0}^t\|\mathbf{e}_n(t)\|_1\mathrm{d}s.
\end{equation}

Let $M=(m+1)Le^{LT}.$ 
From (\ref{vim24}) we obtain successively:

\noindent
For $n=0$ 
$$
\|\mathbf{e}_1(t)\|_1\le M \int_{t_0}^t\|\mathbf{e}_0(s)\|_1\mathrm{d}s\le
 M(t-t_0)\|\mathbf{e}_0\|_{\infty} \ \Rightarrow\ \|\mathbf{e}_1\|_{\infty}\le MT\|\mathbf{e}_0\|_{\infty}.
$$

\noindent
For $n=1$ 
$$
\|\mathbf{e}_2(t)\|_1\le M \int_{t_0}^t\|\mathbf{e}_1(s)\|_1\mathrm{d}s\le
 \frac{M^2 (t-t_0)^2}{2}\|\mathbf{e}_0\|_{\infty}\ \Rightarrow\ \|\mathbf{e}_2\|_{\infty}\le \frac{M^2T^2}{2}
\|\mathbf{e}_0\|_{\infty}.
$$

\noindent
Inductively, it results that
$$
\|\mathbf{e}_n(t)\|_1\le M\int_{t_0}^t\|\mathbf{e}_{n-1}(s)\|_1\mathrm{d}s\le 
\frac{M^n (t-t_0)^n}{n!}\|\mathbf{e}_0\|_{\infty}\ \Rightarrow\ \|\mathbf{e}_n\|_{\infty}\le \frac{M^nT^n}{n!}
\|\mathbf{e}_0\|_{\infty}.
$$
and hence $\lim_{n\rightarrow\infty}\|\mathbf{e}_n\|_{\infty}=0.\quad\rule{5pt}{5pt}$

\section{Computational results}

Although the VIM may be implemented for symbolic computation our experiments are disappointing
for nonlinear equations. As an example the
 \textit{Mathematica} procedure
\scriptsize
\lstset{language=Mathematica}
\begin{lstlisting}
In[1]:=
  VIM[f, U0_, m_] := Module[{V, U = U0, df, Lambda},
  df[t_, x_] := D[f[t, x], x];
  Lambda[U_] := -Exp[
     Integrate[df[w, x] /. {x -> U, t -> w}, {w, s, t}]];
  For[i = 0, i < m, i++,
   V = U + 
     Integrate[Lambda[U] ((D[U, t] - f[t, U]) /. t -> s), {s, 0, t}]; 
   U = V; Clear[V]]; U]
\end{lstlisting}
\normalsize
solves the initial value problem (\ref{vim4})-(\ref{vim5}).

 \begin{exemplu}The problem
$$
\begin{array}{l}
x'(t)=2x(t)+t\\
x(0)=0
\end{array}
$$
\end{exemplu}
with the solution $x(t)=\frac{1}{4}(e^{2t}-2t-1)$
is solved in an iteration,
\scriptsize
\lstset{language=Mathematica}
\begin{lstlisting}
In[1]:= f[t_, x_] := 2 x + t
In[2]:= U0 = 0;
In[3]:= VIM[f,U0,1]
Out[3]= (1/4)*(-1 + E^(2*t) - 2*t)
\end{lstlisting}
\normalsize
but for a nonlinear differential equation it does not give an acceptable result  in a reasonable time.

Much better results we have obtained with numerical computation. The relation (\ref{vim8}) is transformed
into
$$ 
u_{n+1}(t)=\int_{t_0}^t\left(f(s,u_n(s))-f_x(s,u_n(s))u_n(s)\right)e^{\int_s^tf_x(\tau,u_n(\tau))\mathrm{d}\tau}\mathrm{d}s+
$$
$$
+e^{\int_0^tf_x(\tau,u_n(\tau))\mathrm{d}\tau}x^0.
$$ 
$u_n$ will be computed as a first order spline function defined by the points $(t_i,u_n(t_i))_{0\le i\le m}.$ 
$(t_i)_{0\le i\le m}$ is an equidistant grid on $[t_0,t_f].$ The integrals are computed with the trapezoidal rule.


The following \textit{Scilab,} \cite{10}, code solves the initial value problem (\ref{vim4})-(\ref{vim5})
\scriptsize
\lstset{language=Matlab}
\begin{lstlisting}
function [t,u_old,er,iter]=vim(f,df,x0,t0,tf,m,nmi,tol)
    t=linspace(t0,tf,m)
    h=(tf-t0)/(m-1)
    u_old=x0*ones(1,m)
    sw=%t
    iter=0
    while sw do
        iter=iter+1
        u_new=x0*ones(1,m)
        f0=zeros(1,m)
        df0=zeros(1,m)
        for j=1:m do
            f0(j)=f(t(j),u_old(j))
            df0(j)=df(t(j),u_old(j))
        end    
        for i=2:m do
            z=zeros(1,m)
            z(i)=0
            for j=i-1:-1:1 do
                z(j)=0.5*h*(df0(j)+df0(j+1))+z(j+1)
            end
            w=(f0-df0.*u_old).*exp(z)
            s=w(1)+w(i)
            if i>2 then
                for j=2:i-1 do
                    s=s+2*w(j)
                end
            end
            u_new(i)=0.5*h*s+exp(z(1))*x0
        end
        nrm=norm(u_new-u_old,%inf)
        u_old=u_new
        if nrm<tol | iter>=nmi then
            sw=%f
        end
    end 
    if nrm<tol then
        er=0
    else
        er=1
    end 
endfunction
\end{lstlisting}
\normalsize
For the above example
\scriptsize
\lstset{language=Matlab}
\begin{lstlisting}
    deff('y=f(t,x)','y=2*x+t')
    deff('y=df(t,x)','y=2')
    x0=0
    t0=0
    tf=1
    m=100
    nmi=5
    tol=1e-5
\end{lstlisting}
\normalsize
we obtained $\max_{0\le i\le m}|u_2(t_i)-x(t_i)|\approx 0.0000713$ and $\max_{0\le i\le m}|u_2(t_i)-u_1(t_i)|\le tol=10^{-5}.$

\begin{exemplu}
The problem
$$
\begin{array}{l}
x'(t)=x^2(t)+1\\
x(0)=0
\end{array}
$$
\end{exemplu}
has the solution $x(t)=\frac{1-e^{-2t}}{1+e^{-2t}}.$
For $m=100$ the results were  $\max_{0\le i\le m}|u_4(t_i)-x(t_i)|\approx 0.0000293$
and $\max_{0\le i\le m}|u_4(t_i)-u_3(t_i)|\le tol=10^{-5}.$

\end{document}